\documentclass[1p]{elsarticle}
\makeatletter
\gdef\@journal{}
\def\@footnotetwo{}
\makeatother

\usepackage{amsmath, amssymb, amsthm}
\newtheorem{theorem}{Theorem}[section] 
\newtheorem{lemma}[theorem]{Lemma}    

\newtheorem{proposition}[theorem]{Proposition}

\providecommand{\keywords}[1]
{
  \textbf{\textit{Keywords---}} #1
}

\begin{document}
\title
 {Three consecutive near-square squarefree numbers}
\begin{abstract}
In this note, we prove by using T. Estermann's and S. Dimitrov's arguments with an elementary inequality that there are infinitely many $n$ for which all of the numbers $n^2+1,n^2+2$ and $n^2+3$ are squarefree. We also improve the error term slightly in the case of two consecutive numbers of the same form, so that we are able to prove the following asymptotic formula.
\begin{align*}
\sum_{n\le X}\mu^2(n^2+1)\mu^2(n^2+2)\mu^2(n^2+3)\sim\dfrac{7}{18}\prod_{p>3}\left(1-\dfrac{3+\left(\frac{-1}{p}\right)+\left(\frac{-2}{p}\right)+\left(\frac{-3}{p}\right)}{p^2}\right)X.
\end{align*}
\end{abstract}
\author{W. Wongcharoenbhorn}
\ead{w.wongcharoenbhorn@gmail.com}
\address{147 M.1 Baromma-at street
Meaung district, Uttaradit, 53000, Thailand}
\maketitle

\keywords{squarefree number, asymptotic formula, number theory}

\section{Introduction}
We call a positive integer $n$ a \textit{squarefree number} if it is not divisible by a square of a prime. For example, $4$ is not squarefree as $2^2$ divides it, and $6$ is as it is $2\times 3$, while $12$ is not because it is $2^2\times 3$. By convention, we put $1$ to be the smallest squarefree number. Moreover, it is well-known that there are infinitely many squarefree numbers. To be more precise it is proven the following summation formula (see for example \cite{Tenenbaum}).
\begin{align*}
\sum_{n\le X}\mu^2(n)=\dfrac{6X}{\pi^2}+O(\sqrt{X}),
\end{align*}
where $\mu$ is the M\"{o}bius function so that $\mu^2(n)=1$ if $n$ is squarefree and $0$ otherwise.

For simplicity, we may coin the term $\mathcal L\mu$-consecutive formula for the function $f$ (for some positive integer $\mathcal L$) if there exists a formula as follows. 

For some positive value $\sigma\le 1$,
$$\sum_{n\le X}\mu^2(f(n))\mu^2(f(n)+1)\cdots \mu^2(f(n)+\mathcal L-1)  = \sigma X+O(g(X)),$$
where $g(X)=o(X)$, and a non-constant $f(n)$ is some polynomial in $\mathbb Z[n]$.

T. Estermann \cite{Estermann} proved for the first time the infinitude of squarefree numbers of the special form $n^2+1$, i.e. the $1\mu$-consecutive formula, with $g(X)=X^{\frac{2}{3}}\log X$. Recently, in 2021 in \cite{Dimitrov}, S. Dimitrov showed that there are infinitely many pairs $n^2+1, n^2+2$ that are squarefree, and provided that the error term $g(X)$ for $2\mu$-consecutive formula, for $f(n)=n^2+1$, can be $O(X^{\frac{8}{9}+\varepsilon})$. We employ the arguments from T. Estermann's and S. Dimitrov's work to prove the result stated in the abstract of this note.

On the direction of $1\mu$-consecutive formula for $f(n)=n^2+1$, Heath-Brown \cite{Heath-Brown} improved the error estimate to $O(X^{\frac{7}{12}+\varepsilon})$. We are not interested in this direction, but we improve the bound slightly after the work of S. Dimitrov. Finally, it is worth noting that there are no four consecutive squarefree numbers of any form, since one of any four consecutive numbers is divisible by $2^2$. Therefore, three is the longest length for consecutive squarefree numbers of any type and hence, $\mathcal L\le 3$.

We will prove here the exact $3\mu$-consecutive formula for $f(n)=n^2+1$, in our main proposition and its weaker form by an elementary inequality. In particular, we prove
\begin{proposition}
We have that
\label{prop:main}
\begin{align*}
\sum_{n\le X}&\mu^2(n^2+1)\mu^2(n^2+2)\mu^2(n^2+3)\\
&=\dfrac{7}{18}\prod_{p>3}\left(1-\dfrac{3+\left(\frac{-1}{p}\right)+\left(\frac{-2}{p}\right)+\left(\frac{-3}{p}\right)}{p^2}\right)X+O(X^{\frac{3}{4}+\varepsilon}),
\end{align*}
as $X$ tends to infinity.
\end{proposition}

In summary, the proposition is obtained by combining the summation formula for $1\mu$- and $2\mu$-consecutive for $f(n)=n^2+1$ and the arguments used in T.Estermann's work to reduce the error term. It is clear that the corollary of this proposition (and its weaker form) is the infinitude of three consecutive squarefree numbers of the form $n^2+1,n^2+2$, and $n^2+3$. 

The organization of this note is as follows. Firstly, in section \ref{sec:prelim} we introduce some notations and terms used. Then, in section \ref{sec:proplem} we present important propositions and lemmas. Also, we clarify here what do we mean by the weaker form of the main proposition \ref{prop:main} and show how we use the elementary inequality to imply it. Finally, in section \ref{sec:main} we provide proofs of the propositions stated in section \ref{sec:proplem}, which are mainly the summation formula for $2\mu$-consective for $f(n)=n^2+1$, and the completion of the proof of our main proposition \ref{prop:main}.

\section{Preliminaries}
\label{sec:prelim}

Let us denote $\Gamma(X)=\sum_{n\le X}\mu^2(n^2+1)\mu^2(n^2+2)\mu^2(n^2+3)$. For $m,k,j\in\left\{1,2,3\right\}$ with $k<j$, we denote constants
\begin{align*}
\Phi_k &= c_k\prod_{p>3}\left(1-\dfrac{1+\left(\frac{-k}{p}\right)}{p^2}\right),\\
\Theta_m &= c_{k,j}\prod_{p>3}\left(1-\dfrac{2+\left(\frac{-k}{p}\right)+\left(\frac{-j}{p}\right)}{p^2}\right),
\end{align*}
where $c_1=1, c_2=7/9$ and $c_3=1/2$, $c_{1,2}=8/9, c_{1,3}=1/2$ and $c_{2,3}=7/18$, and $\left(\frac{\cdot}{p}\right)$ is the Legendre symbol over an odd prime $p$,.

Here, the constants $\Theta_1,\Theta_2,\Theta_3$ correspond to the tuple $(k,j)=(1,2),(1,3),(2,3)$ respectively. Important sets involving in this note for $\phi=1,2,3$ are denoted by
\begin{align*}
N'_\phi(d) = \left\{1\le n\le d^2: n^2+\phi\equiv 0\pmod{d^2}\right\}.
\end{align*}

Furthermore, though we need no use of exponential sum, we let $B_1(t)=\left\langle\, t\right\rangle -1/2$, where $\left\langle\, t\right\rangle$ is the fractional part of $t$, and as a convention we write $t=\lfloor t\rfloor+\left\langle\,t\right\rangle $, where $\lfloor t\rfloor$ is the largest integer not exceeding $t$.

For the arithmetic functions $\mu,\tau$ we mean the M\"{o}bius function and the divisor function respectively. As a convention in number theory, we mean by $n,m,k,j$ the integers usually positive, $p$ denotes a prime, and $X$ is always a large positive number, while $\varepsilon$ is small positive number and can be varied throughout the note.

\section{Propositions and lemmas}
\label{sec:proplem}
In this section, we address some propositions and lemmas used in this essay. The $\Phi$'s and $\Theta$'s appearing in the following two propositions are as defined in section \ref{sec:prelim}.
\begin{proposition}
\label{prop:sum_k}
For $k\in\left\{1,2,3\right\}$,
\begin{align*}
\sum_{n\le X}\mu^2(n^2+k) = \Phi_kX+O\Big(X^{\frac{2}{3}}\log X\Big).
\end{align*}
\end{proposition}

\begin{proof}
This is the partial result of \cite{Estermann}.
\end{proof}

Here, we present another main proposition that the proof will be deferred to the next section. We improve the bound on the error term slightly as follows.
\begin{proposition}
\label{prop:sum_kj}
For $m,k,j\in\left\{1,2,3\right\}$ with $k<j$,
\begin{align*}
\sum_{n\le X}\mu^2(n^2+k)\mu^2(n^2+j)= \Theta_m X+O\Big(X^{\frac{2}{3}+\varepsilon}\Big).
\end{align*}
\end{proposition}

\begin{lemma}
\label{lem:sum_B1}
If we let $$\mathcal{B}(D)=\sum_{D\le d<2D}\sum_{n\in N'_{\phi}(d)} \left( B_1\left(\dfrac{-n}{d^2}\right)-B_1\left(\dfrac{X-n}{d^2}\right)\right),$$
then for $\delta>0$, and for any integer $n\le X$, we obtain that for any small $\varepsilon>0$
\begin{align*}
\mathcal{B}(D)\ll
\begin{cases}
X^{2\delta+\varepsilon}, \text{ if $D>X^{1-\delta}$},\\
X^{1-\delta+\varepsilon}, \text{ otherwise}.
\end{cases}
\end{align*}
\end{lemma}
\begin{proof}
Here, we denote $\phi=1,2,3$. For $D>X^{1-\delta}$, if we have $n^2-kd^2=-\phi$, then $1\le k\le X^{2\delta}$ as $n\le X$ and $d>X^{1-\delta}$. The left-hand side is bounded above by
\begin{align*}
\sum_{D\le d<2D}\sum_{n\in N'_{\phi}(d)} 1 \ll \sum_{k\le X^{2\delta}}\#\left\{(n,d):n\le X, n^2-kd^2=-\phi\right\}:=\sum_{k\le X^{2\delta}} A(k),
\end{align*}

where we denote the number of solutions to $n^2-kd^2=-\phi$ for fixed $k$ be $A(k)$. Then, we may argue as in \cite{Estermann} that $A(k)\ll X^{\varepsilon}$ since $kd^2\le X^2$, and hence the first assertion of the lemma follows.

For $D\le X^{1-\delta}$, we again have that $\mathcal{B}(D)$ is bounded by
\begin{align*}
\sum_{D\le d<2D}\sum_{n\in N'_{\phi}(d)} 1 \ll \sum_{d\le X^{1-\delta}} \#\left\{n\le d^2: n^2+\phi\equiv 0\pmod{d^2}\right\}\ll \sum_{d\le X^{1-\delta}}2^{\omega(d)},
\end{align*}
since the number of solutions to the quadratic polynomial is at most $2^{\omega(d^2)}=2^{\omega(d)}$, where we denote $\omega(k)$ be the number of different prime factors of $k$. 

We have the formula that for any $d$, $2^{\omega(d)}=\sum_{k|d}\mu^2(k)$. To see this, we first note that $2^{\omega}$ is multiplicative, i.e., for $d_1,d_2$ such that $\gcd(d_1,d_2)=1$, $2^{\omega(d_1d_2)}=2^{\omega(d_1)}\cdot 2^{\omega(d_2)}$, and that $\sum_{k|d}\mu^2(k)$ is also multiplicative. Thus, to complete the proof of the formula, it is sufficient to show for a prime power $p^{\nu}$ for $\nu\ge 1$ that
$$2^{\omega(p^{\nu})}=2=\mu^2(1)+\mu^2(p)=\sum_{k|p^{\nu}}\mu^2(k).$$

Therefore, by the above formula we have that
\begin{align*}
\sum_{d\le Y} 2^{\omega(d)}=\sum_{d\le Y}\sum_{k|d}\mu^2(k)\ll\sum_{k\le Y}\mu^2(k)\cdot\left(\dfrac{Y}{k}\right) \ll Y\log Y,
\end{align*}
so that the lemma follows.
\end{proof}

\begin{lemma}
\label{lem:compute-const}
For $\Phi_1,\Phi_2,\Phi_3$ and $\Theta_1,\Theta_2,\Theta_3$ denoted as in the preliminaries section, then
\begin{align*}
\Phi_1&\in(0.889617,0.894956), \Phi_2\in(0.746484,0.750964), \Phi_3\in(0.464691,0.467481),\\
\Theta_1&\in(0.668036,0.672046), \Theta_2\in(0.415817,0.418313), \Theta_3\in(0.348955,0.351050).
\end{align*}.
\end{lemma}

\begin{proof}
This can be done by the computation and the following inequalities.
\begin{align*}
\prod_{p> 1000}\left(1-\dfrac{4}{p^2}\right)\prod_{3<p\le 1000}\left(1-\dfrac{\rho(p)}{p^2}\right)\le\prod_{p>3}\left(1-\dfrac{\rho(p)}{p^2}\right)<\prod_{3<p\le 1000}\left(1-\dfrac{\rho(p)}{p^2}\right),
\end{align*}
where $|\rho(p)|\le 4$ denotes the number in the numerator of each product. 

Now we determine the number
\begin{align*}
\prod_{p> 1000}\left(1-\dfrac{4}{p^2}\right)\ge \prod_{n>1000}\left(1-\dfrac{4}{n^2}\right)=\dfrac{1}{6}\prod_{3\le n\le 1000}\left(1-\dfrac{4}{n^2}\right)^{-1}>1/1.006,
\end{align*}
by computation. Note that, the product of terms involving primes not greater than $1000$ can be computed, and that every $\Phi_m$'s and $\Theta_k$'s satisfy the property for each numerator in each product. For if we denote the numerator be $\rho(p)$, we have $|\rho(p)|\le 4$. Hence, we are done.
\end{proof}

Then, we prove here that the propositions and lemmas above (in this section) imply the weaker form of proposition \ref{prop:main}, which asserts an existence of a positive constant $\alpha>0.1477$, such that $\Gamma(X)>\alpha X$ for all large $X$. This can be done without the reduction of the error term in the case of two consecutive numbers. 

Firstly, we obtain by the Power mean inequality
\begin{align*}
&\left(\dfrac{1}{\lfloor X\rfloor}\cdot\sum_{n\le X}\left(\mu^2(n^2+1)+\mu^2(n^2+2)\mu^2(n^2+3)\right)^2\right)^{\frac{1}{2}}\\
&\ge \left(\dfrac{1}{\lfloor X\rfloor}\cdot\sum_{n\le X}\left(\mu^2(n^2+1)+\mu^2(n^2+2)\mu^2(n^2+3)\right)\right).
\end{align*}

By expanding the above inequality, using proposition \ref{prop:sum_k}, proposition \ref{prop:sum_kj} (where we may be allowed to use the error term to be $o(X)$), and lemma \ref{lem:compute-const}, we obtain the lower bound
\begin{align*}
\Gamma(X)\ge\frac{1}{2}\left(\Big(\Phi_1+\Theta_{3}\Big)^{2}-(\Phi_1+\Theta_3)\right)X+o(X)> \alpha X,
\end{align*}
where the inequality holds for all sufficiently large $X$.

Let us define some terms further. For $k<j$ and $k,j\in\left\{1,2,3\right\}$ we denote (sometimes we omit the subscript $k,j$ for brevity)
\begin{align*}
\Gamma_{k,j}&=\sum_{n\le X}\mu^2(n^2+k)\mu^2(n^2+j),\\
S_{k,j}(q_1,q_2) &= \left\{n\in\mathbb{N}: 1\le n\le q_1q_2, n^2+k\equiv 0\pmod{q_1}, n^2+j\equiv 0\pmod{q_2}\right\},\\
\lambda_{k,j}(q_1,q_2) &= \sum_{n\in S_{k,j}(q_1,q_2)}1.
\end{align*}

\begin{lemma}
\label{lem:multiplicative}
For $k,j\in\left\{1,2,3\right\}$ with $k<j$, denote $\lambda_{k,j}=\lambda$. We have for $\gcd(q_1q_2,q_3q_4)=\gcd(q_1,q_2)=\gcd(q_3,q_4)=1,$
$$\lambda(q_1q_2,q_3q_4)=\lambda(q_1,q_3)\lambda(q_2,q_4).$$

\begin{proof}
We first note by the Chinese remainder theorem that $n$ on the left-hand side of the summation below is from the combination of all possibilities of the following congruences.
\begin{align*}
n&\equiv n_{1,a_1}\pmod{q_1},
n\equiv n_{2,a_2}\pmod{q_2},
n\equiv n_{3,a_3}\pmod{q_3},
n\equiv n_{4,a_4}\pmod{q_4},
\end{align*}
where $n_{1,a_1}, n_{2,a_2}, n_{3,a_3}$ and $n_{4,a_4}$ satisfy
\begin{align*}
n_{1,a_1}^2+k&\equiv 0\pmod{q_1},\\
n_{2,a_2}^2+k&\equiv 0\pmod{q_2},\\
n_{3,a_3}^2+j&\equiv 0\pmod{q_3},\\
n_{4,a_4}^2+j&\equiv 0\pmod{q_4}.
\end{align*}

Thus, we can apply the multiplication law of counting to obtain
$$\sum_{\substack{n\le q_1q_2q_3q_4\\n^2+k\equiv 0 (q_1q_2)\\n^2+j\equiv 0 (q_3q_4)}}1=\left(\sum_{\substack{n\le q_1\\n^2+k\equiv 0(q_1)}}1\right)\left(\sum_{\substack{n\le q_2\\n^2+k\equiv 0(q_2)}}1\right)\left(\sum_{\substack{n\le q_3\\n^2+j\equiv 0(q_3)}}1\right)\left(\sum_{\substack{n\le q_4\\n^2+j\equiv 0(q_4)}}1\right).$$
Then, we rearrange the product on the right-hand side above and employ again the multiplication law. This completes the proof of the lemma.
\end{proof}
\end{lemma}

\section{Main result}
\label{sec:main}

In this section, we now prove proposition \ref{prop:sum_kj} and our main proposition \ref{prop:main}. Using almost the same argument as in \cite{Dimitrov}, we obtain
\begin{align*}
\Gamma_{k,j}(X) &= \sum_{\substack{d_1,d_2\\ \gcd(d_1,d_2)=1}}\mu(d_1)\mu(d_2)\sum_{\substack{n\le X\\ n^2+k\equiv 0(d_1^2)\\ n^2+j\equiv 0(d_2^2)}} 1=\Gamma_1(X)+\Gamma_2(X),
\end{align*}
where $\sqrt{X}\le z < X$ will be chosen later, and
\begin{align*}
\Gamma_1(X) &= \sum_{\substack{d_1d_2\le z\\ \gcd(d_1,d_2)=1}} \mu(d_1)\mu(d_2)\Sigma(X,d_1^2,d_2^2),\\
\Gamma_2(X) &= \sum_{\substack{d_1d_2> z\\ \gcd(d_1,d_2)=1}} \mu(d_1)\mu(d_2)\Sigma(X,d_1^2,d_2^2),\\
\Sigma(X,d_1^2,d_2^2) &= \sum_{\substack{n\le X\\ n^2+k\equiv 0(d_1^2)\\ n^2+j\equiv 0(d_2^2)}}1.
\end{align*}

\subsection{Estimation of $\Gamma_{1}(X)$}
Suppose that $q_1=d_1^2, q_2=d_2^2$ where $d_1,d_2$ are squarefree numbers with $\gcd(q_1,q_2)=1$ and $d_1d_2\le z$. Denote
$$\Omega(X,q_1,q_2,n)=\sum_{\substack{m\le X\\ m\equiv n(q_1q_2)}}1.$$

We obtain further that
\begin{align*}
\Sigma(X,q_1,q_2) &= \sum_{n\in S(q_1,q_2)}\Omega(X,q_1,q_2,n),\\
\Omega(X,q_1,q_2,n) &= \dfrac{X}{q_1q_2}+O(1).
\end{align*}
Therefore, we have
\begin{align*}
\Sigma(X,q_1,q_2) = X\cdot\dfrac{\lambda(q_1,q_2)}{q_1q_2}+O(\lambda(q_1,q_2)).
\end{align*}
Furthermore, note that $\lambda(q_1,q_2)\ll\tau(q_1q_2)\ll X^{\varepsilon}$. We see that
\begin{align*}
\sum_{\substack{d_1d_2\le z\\ \gcd(d_1,d_2)=1}} 1\ll \sum_{n\le z}\tau(n)\ll z\log z,
\end{align*}
and that $\log z$ is absorbed into $X^{\varepsilon}$ since $z\le X$. Now, we obtain by plugging in the above equation and inequality
\begin{align*}
\Gamma_1(X) &= X\left(\sum_{\substack{d_1d_2\le z\\ \gcd(d_1,d_2)=1}}\dfrac{\mu(d_1)\mu(d_2)\lambda(d_1^2,d_2^2)}{d_1^2d_2^2}\right)+O(zX^{\varepsilon})\\
&= \sigma X-X\left(\sum_{\substack{d_1d_2>z\\ \gcd(d_1,d_2)=1}}\dfrac{\mu(d_1)\mu(d_2)\lambda(d_1^2,d_2^2)}{d_1^2d_2^2}\right)+O(zX^{\varepsilon}),
\end{align*}
where we denote
\begin{align*}
\sigma =\sum_{\substack{d_1,d_2\ge 1\\ \gcd(d_1,d_2)=1}}\dfrac{\mu(d_1)\mu(d_2)\lambda(d_1^2,d_2^2)}{d_1^2d_2^2}.
\end{align*}

Now we determine the error term occurs in the second term attached to $X$ above. Note that we have (where the last inequality can be found in \cite{Apostol})
\begin{align*}
\sum_{\substack{d_1d_2>z\\ \gcd(d_1,d_2)=1}}\dfrac{\mu(d_1)\mu(d_2)\lambda(d_1^2,d_2^2)}{d_1^2d_2^2}\ll \sum_{\substack{d_1d_2>z\\ \gcd(d_1,d_2)=1}}\dfrac{(d_1d_2)^{\varepsilon}}{d_1^2d_2^2}\ll \sum_{\substack{n>z}} \dfrac{1}{n^{2-\varepsilon}}\ll \dfrac{1}{z^{1-\varepsilon}}.
\end{align*}

Whence, for a positive $\varepsilon<1$, we need only to show that $\sigma$ corresponds to the constant in proposition \ref{prop:sum_kj}. By lemma \ref{lem:multiplicative}, we have $\lambda(d_1^2,d_2^2)=\lambda(d_1^2,1)\lambda(1,d_2^2)$. Then, by letting
\begin{align*}
f(d_1,d_2) = 
\begin{cases}
1 & \text{if $\gcd(d_1,d_2)=1$},\\
0 & \text{if $\gcd(d_1,d_2)>1$},
\end{cases}
\end{align*}
we have that, as both sums over $d_1$ and $d_2$ are absolute convergent,
\begin{align*}
\sigma &= \sum_{d_1\ge 1}\dfrac{\mu(d_1)\lambda(d_1^2,1)}{d_1^2}\sum_{d_2\ge 1}\dfrac{\mu(d_2)\lambda(1,d_2^2)}{d_2^2}f(d_1,d_2)\\
&= \prod_p\left(1-\dfrac{\lambda(1,p^2)}{p^2}\right)\sum_{d_1\ge 1}\dfrac{\mu(d_1)\lambda(d_1^2,1)}{d_1^2}\prod_{p|d_1}\left(1-\dfrac{\lambda(1,p^2)}{p^2}\right)^{-1}\\
&= \prod_{p}\left(1-\dfrac{\lambda(1,p^2)}{p^2}\right)\prod_{p}\left(1-\dfrac{\lambda(p^2,1)}{p^2}\left(1-\dfrac{\lambda(1,p^2)}{p^2}\right)^{-1}\right)\\
&= \prod_{p}\left(1-\dfrac{\lambda(p^2,1)+\lambda(1,p^2)}{p^2}\right)=c_{k,j}\prod_{p>3}\left(1-\dfrac{2+\left(\frac{-k}{p^2}\right)+\left(\frac{-j}{p^2}\right)}{p^2}\right).
\end{align*}
Therefore, we are left to bound the error term to be $O\Big(X^{\frac{2}{3}+\varepsilon}\Big)$.

\subsection{Estimation of $\Gamma_{2}(X)$}
By splitting into dyadic ranges as in \cite{Dimitrov}, we obtain

\begin{align*}
\Gamma_2(X)\ll \log^2{X}\sum_{n\le X}\sum_{\substack{D\le d<2D\\ n^2+k\equiv 0(d^2)}}\sum_{\substack{D'\le d'<2D'\\ n^2+j\equiv 0(d'^2)}}1,
\end{align*}
where $D, D'$ satisfy the conditions
$$\dfrac{1}{2}\le D,D'\le \sqrt{X^2+\max(k,j)}, DD'>\dfrac{z}{4}.$$
Then for $\phi=1,2,3$, we obtain by noting that the number of solutions to $n^2+\phi=kd^2$ for $d\in[D,2D)$ and fixed $n\le X$ is at most $\tau(n^2+\phi)\ll n^\varepsilon \ll X^{\varepsilon}$,
\begin{align*}
\Gamma_2(X)\ll X^{\varepsilon}\Big(\Sigma_k+\Sigma_j\Big),
\end{align*}
where we may denote $D,D'$ interchangeably that corresponds to $\phi$ and here
\begin{align*}
\Sigma_\phi =\sum_{n\le X}\sum_{\substack{D\le d<2D\\ n^2+\phi\equiv 0(d^2)}}1.
\end{align*}

Then, we estimate the above term for fixed $\phi$ as follows. Remark that the condition $n\le X$ always holds.
\begin{align*}
\Sigma_\phi &= \sum_{D\le d<2D}\sum_{n\in N'_{\phi}(d)}\sum_{\substack{m\le X\\ m\equiv n(d^2)}}1 = \sum_{D\le d<2D}\sum_{n\in N'_{\phi}(d)}\left(\left\lfloor \dfrac{X-n}{d^2} \right\rfloor-\left\lfloor \dfrac{-n}{d^2} \right\rfloor\right)\\
&= \sum_{D\le d<2D}\sum_{n\in N'_{\phi}(d)} \left(\dfrac{X}{d^2}+B_1\left(\dfrac{-n}{d^2}\right)-B_1\left(\dfrac{X-n}{d^2}\right)\right).
\end{align*}

Hence, we have
\begin{align*}
\Sigma_{\phi}\ll X^{1+\varepsilon}D^{-1}+\left|\Sigma^{(0)}_{\phi}\right|,
\end{align*}
where we denote
\begin{align*}
\Sigma^{(0)}_{\phi} &= \sum_{D\le d<2D}\sum_{n\in N'_{\phi}(d)} \left(B_1\left(\dfrac{-n}{d^2}\right)-B_1\left(\dfrac{X-n}{d^2}\right)\right).
\end{align*}

By lemma \ref{lem:sum_B1}, we obtain
\begin{align*}
\Sigma_{\phi}^{(0)} \ll X^{2\delta+\varepsilon}+X^{1-\delta+\varepsilon}.
\end{align*}

Therefore, after replacing the constraint we have that the error term from both terms $\Gamma_1(X),\Gamma_2(X)$ is bounded above by,
\begin{align*}
\ll 
zX^{\varepsilon}+X^{1+\varepsilon}z^{-\frac{1}{2}}+X^{1-\delta+\varepsilon}+X^{2\delta+\varepsilon}.
\end{align*}

By choosing the optimal choice $z=X^{\frac{2}{3}}$ and $\delta=\frac{1}{3}$, we obtain the error term in propostition \ref{prop:sum_kj} to be $O(X^{\frac{2}{3}+\varepsilon})$. Hence, we acquire proposition \ref{prop:sum_kj}, thereby obtaining the weak form of proposition \ref{prop:main}. It is worth noting that this computation can be extended to any polynomial of the form $n^2+\ell$. for $\ell=1,2\dots$ be any positive integer.

\subsection{The asymptotic formula of $\Gamma(X)$}

Now, we turn ourselves to the proof of proposition \ref{prop:main}. Firstly, we denote the conditions used in the summation as follows.
\begin{align*}
&(1) \gcd(d_1,d_2)=\gcd(d_2,d_3)=\gcd(d_3,d_1)=1,\\
&(2) n^2+i\equiv 0\pmod{d_i^2}, \forall i=1,2,3
\end{align*}
Using the arguments as in the previous case, we obtain that
\begin{align*}
\Gamma(X) &= \sum_{\substack{d_1,d_2,d_3\\ (1)}}\mu(d_1)\mu(d_2)\mu(d_3)\sum_{\substack{n\le X\\ (2)}} 1:=\Gamma_1(X)+\Gamma_2(X),
\end{align*}
where $z\le X$ will be chosen later, and
\begin{align*}
\Gamma_1(X) &= \sum_{\substack{d_1d_2d_3\le z\\ (1)}} \mu(d_1)\mu(d_2)\mu(d_3)\Sigma(X,d_1^2,d_2^2,d_3^2),\\
\Gamma_2(X) &= \sum_{\substack{d_1d_2d_3> z\\ (1)}} \mu(d_1)\mu(d_2)\mu(d_3)\Sigma(X,d_1^2,d_2^2,d_3^2),\\
\Sigma(X,d_1^2,d_2^2,d_3^2) &= \sum_{\substack{n\le X\\ (2)}}1.
\end{align*}

Then, we also coin the term analogous to the previous $\lambda_{k,j}$'s that are used in the case of two consecutive terms,
\begin{align*}
\lambda(q_1,q_2,q_3)=\sum_{\substack{n\le q_1q_2q_3\\ (2)}}1,
\end{align*}
where $q_i=d_i^2,\forall i\in\left\{1,2,3\right\}$. So that we have $\lambda(q_1,q_2,q_3)=\lambda(q_1,1,1)\lambda(1,q_2,1)\lambda(1,1,q_3)$. This can be proven in the same way as in the proof of lemma \ref{lem:multiplicative}.

Therefore, we obtain similarly that
\begin{align*}
\Sigma(X,q_1,q_2,q_3)=\dfrac{X\cdot\lambda(q_1,q_2,q_3)}{q_1q_2q_3}+O(\lambda(q_1,q_2,q_3)),
\end{align*}
and that (since $\lambda(q_1,q_2,q_3)\ll X^{\varepsilon}$)
\begin{align*}
\Gamma_1(X)=\sigma X+O(zX^{\varepsilon}),
\end{align*}
where $\sigma$ here denotes
\begin{align*}
\sigma = \sum_{d_1\ge 1}\dfrac{\mu(d_1)\lambda(d_1^2,1,1)}{d_1^2}\sum_{d_2\ge 1}\dfrac{\mu(d_2)\lambda(1,d_2^2,1)}{d_2^2}\sum_{d_3\ge 1}\dfrac{\mu(d_3)\lambda(1,1,d_3^2)}{d_3^2}\tilde{f}(d_1,d_2,d_3),
\end{align*}
for which $\tilde{f}(a,b,c)=f(a,b)f(b,c)f(c,a)$, i.e., $\tilde{f}(a,b,c)=1$ if $a,b,c$ are pairwise coprimes and $0$ otherwise. We also have further that $\sigma$ equals
\begin{align*}
&\sum_{d_1\ge 1}\dfrac{\mu(d_1)\lambda(d_1^2,1,1)}{d_1^2}\sum_{d_2\ge 1}\dfrac{\mu(d_2)\lambda(1,d_2^2,1)}{d_2^2}\sum_{d_3\ge 1}\dfrac{\mu(d_3)\lambda(1,1,d_3^2)}{d_3^2}\tilde{f}(d_1,d_2,d_3)\\
&= \prod_p\left(1-\dfrac{\lambda(1,p^2,1)+\lambda(1,1,p^2)}{p^2}\right)\sum_{d_1\ge 1}\dfrac{\mu(d_1)\lambda(d_1^2,1,1)}{d_1^2}\prod_{p|d_1}\left(1-\dfrac{\lambda(1,p^2,1)+\lambda(1,1,p^2)}{p^2}\right)^{-1}\\
&= \prod_p\left(1-\dfrac{\lambda(p^2,1,1)+\lambda(1,p^2,1)+\lambda(1,1,p^2)}{p^2}\right)\\
&=\dfrac{7}{18}\prod_{p>3}\left(1-\dfrac{3+\left(\frac{-1}{p}\right)+\left(\frac{-2}{p}\right)+\left(\frac{-3}{p}\right)}{p^2}\right).
\end{align*}

Then, we are left to bound $\Gamma_2(X)$. For the conditions $$\dfrac{1}{2}\le D_1,D_2,D_3\le \sqrt[3]{X^2+3}, D_1D_2D_3>\dfrac{z}{8},$$ we have (by splitting into dyadic ranges in three dimensions) that
\begin{align*}
\Gamma_2(X)\ll \log^3 X\sum_{n\le X}\sum_{\substack{D_1\le d_1<2D_1\\ n^2+1\equiv 0(d_1^2)}}\sum_{\substack{D_2\le d_2<2D_2\\ n^2+2\equiv 0(d_2^2)}}\sum_{\substack{D_3\le d_3<2D_3\\ n^2+3\equiv 0(d_1^3)}}1.
\end{align*}
Thus, with the same arguments as in the case of two consecutive terms, we also obtain the following inequality.
$$\Gamma_2(X)\ll X^{\varepsilon}\Big(\Sigma_1+\Sigma_2+\Sigma_3\Big),$$
where now we denote for $\phi=1,2,3$ corresponding to $D_{\phi}=D$,
\begin{align*}
\Sigma_\phi =\sum_{n\le X}\sum_{\substack{D\le d<2D\\ n^2+\phi\equiv 0(d^2)}}1.
\end{align*}
This is bounded by $X^{\frac{2}{3}+\varepsilon}$, by lemma \ref{lem:sum_B1}, for any small $\varepsilon>0$. The error term combined for the case of three consecutive terms thus is
$$\ll zX^{\varepsilon}+X^{1+\varepsilon}z^{-\frac{1}{3}}+X^{\frac{2}{3}+\varepsilon}.$$

Now, the optimal bound occurs when $z=X^{\frac{3}{4}}$ for which we obtain the error estimate in our main proposition \ref{prop:main}. Hence, the proof of our main proposition is completed.

It is worth noting that a similar argument for four or more consecutive terms is impossible. This can be seen as we will have the constant for $\Gamma_1(X)$ of the case of four terms be
\begin{align*}
\prod_{p}\left(1-\dfrac{\lambda(p^2,1,1,1)+\lambda(1,p^2,1,1)+\lambda(1,1,p^2,1)+\lambda(1,1,1,p^2)}{p^2}\right)=0.
\end{align*}

We have the above term as it is possible to show that $d_1^2,d_2^2,d_3^2,d_4^2$ are pairwise coprimes if each $d_i$ corresponds to $n^2+i\equiv 0\pmod{d_i^2},\forall i=1,2,3,4$ with fixed $n$. Here, $\lambda$ is defined in a similar manner as in the case of three consecutive terms.  The above equality holds as for $p=2$, we have $\lambda(p^2,1,1,1)+\lambda(1,p^2,1,1)+\lambda(1,1,p^2,1)+\lambda(1,1,1,p^2)=p^2$. 

In general, for more than four consecutive terms, we have that the main term $\sigma$ (analogous to the cases of two and three consecutive terms) equals to $0$. This can be seen again that $d^2_1,d^2_2,d^2_3,d^2_4$ are pairwise coprimes and that the above product is contained in $\sigma$. Whence, we are unable to conclude the infinitude of the case of four or more consecutive terms of this type by this method.

This research did not receive any specific grant from funding agencies in the public, commercial, or not-for-profit sectors.

\end{document}